\documentclass[11pt]{amsart}

\usepackage{graphicx}
\usepackage{mathptmx}

\usepackage{amscd}
\usepackage{amsxtra}
\usepackage{a4wide}
\usepackage{latexsym}
\usepackage{amssymb}
\usepackage{amsfonts}
\usepackage{amsmath}
\usepackage{amsrefs}
\usepackage{mathrsfs}
\usepackage{upref}
\usepackage{txfonts}
\usepackage{color}
\usepackage{tcolorbox}
\usepackage[utf8]{inputenc}

\usepackage[bookmarksnumbered, colorlinks, plainpages]{hyperref}
\hypersetup{colorlinks=true,linkcolor=red, anchorcolor=green,
	citecolor=cyan, urlcolor=red, filecolor=magenta, pdftoolbar=true}

\usepackage[left=2cm, right=2cm, top=2cm, bottom=2cm]{geometry}

\allowdisplaybreaks

\theoremstyle{plain}
\newtheorem{thm}{Theorem}[section]
\newtheorem{lem}[thm]{Lemma}

\theoremstyle{definition}
\newtheorem{rem}[thm]{Remark}

\newtheorem{defi}[thm]{Definition}

\numberwithin{thm}{section}
\numberwithin{equation}{section}

\newcommand{\zstroke}{%
	\text{\ooalign{\hidewidth -\kern-.4em-\hidewidth\cr$\mathcal Z$\cr}}%
}

\def\supp{\operatorname{supp}}

\def\esup{\operatornamewithlimits{ess\,sup}}

\def\Z{\mathbb Z}

\def\card{\operatorname{card}}

\def\mp{{\mathfrak M}}

\def\I{(0,\infty)}

\begin{document}

\title[{Corrigendum to [Positivity, 22 (1) (2018), 275-299]}]{Corrigendum to "On weighted iterated Hardy-type \\ inequalities" [Positivity, 22 (1) (2018), 275-299]}

\author[R.Ch. Mustafayev]{RZA MUSTAFAYEV}
\address{RZA MUSTAFAYEV, Department of Mathematics, Kamil \"{O}zda\u{g} Faculty of Science, Karamano\u{g}lu Mehmetbey University, 70200, Karaman, Turkey}
\email{rzamustafayev@gmail.com}



\begin{abstract}
We correct a mistake in the paper ["On weighted iterated Hardy-type inequalities", Positivity, 22 (1) (2018), 275-299].
\end{abstract}

\maketitle

\section{Introduction}\label{in}

This is a corrigendum to the paper ["On weighted iterated Hardy-type inequalities", Positivity, 22 (1) (2018), 275-299], which will be referred to as \cite{mus.2017} hereafter. After the paper was published, an error was found in the statement and proof of \cite[Lemma 4.8]{mus.2017}, which caused \cite[Lemma 4.9]{mus.2017} and \cite[Theorem 1.1]{mus.2017} to be erroneous as well. Here we will give a corrected version of \cite[Theorem 1.1]{mus.2017} and its proof. 

We will use the notation and terminology in \cite{mus.2017}. Recall that \cite{mus.2017} deals with the characterization of the inequality
\begin{equation}\label{main}
\bigg( \int_0^{\infty} \bigg( \int_x^{\infty} \bigg( \int_t^{\infty} h \bigg)^q w(t)\,dt
\bigg)^{r / q} u(x)\,dx \bigg)^{1/r}\leq C \,\int_0^{\infty} h v, \quad h \in {\mathfrak M}^+(0,\infty),
\end{equation}
where $0 < q ,\, r < \infty$ and $u,\,v,\,w$ are weight functions on $(0,\infty)$.

\cite[Theorem 1.1]{mus.2017} must be replaced with the following statement.
\begin{thm}\label{main1}
	Let $0 < q,\, r < \infty$ and  $u,\,v,\,w \in {\mathcal W}\I$.
	
	{\rm (a)} Let $1 \le \min\{q,\,r\}$. Then inequality \eqref{main} holds if and only if $F_1 < \infty$, where
	\begin{align*}
		F_1 : & = \esup_{t\in (0,\infty)} \bigg(\int_{0}^{t} u(s) \bigg(\int_{s}^{t} w \bigg)^{{r} / {q}} \, ds\bigg)^{1 / r} \bigg( \esup_{y\in [t,\infty)} v(y)^{-1} \bigg).
	\end{align*}	
	
	Moreover, if $C$ is the best constant in \eqref{main}, then $C \approx F_1$.
	
	{\rm (b)} Let $r < 1 \le q$. Then inequality \eqref{main} holds if and only if $F_2 < \infty$ and $F_3 < \infty$, where
	\begin{align*}
		F_2 & := \bigg(\int_{0}^{\infty}\bigg(\int_{0}^{t} u \bigg)^{r'} u(t) \bigg[\esup_{x\in (t,\infty)} \bigg(\int_{t}^{x} w \bigg)^{{r'} / {q}} \bigg( \esup_{y\in [x,\infty)} v(y)^{-1} \bigg)^{r'} \bigg] \, dt\bigg)^{{1} / {r'}} \\
		F_3 & := \bigg(\int_{0}^{\infty}\bigg(\int_{0}^{t} u(s) \bigg(\int_{s}^{t} w \bigg)^{{r} / {q}} \, ds \bigg)^{r'} u(t) \bigg( \esup_{x\in (t,\infty)} \bigg(\int_{t}^{x} w \bigg)^{{r} / {q}} \bigg( \esup_{y\in [x,\infty)} v(y)^{-1} \bigg)^{r'} \bigg) \, dt \bigg)^{{1} / {r'}}.	
	\end{align*}
	
	Moreover, if $C$ is the best constant in \eqref{main}, then $C \approx F_2 + F_3$.

	{\rm (c)} Let $q < 1 \le r$. Then inequality \eqref{main} holds if and only if $F_1 < \infty$ and $F_4 < \infty$, where
    \begin{align*}
    F_4 & := \esup_{t\in (0,\infty)} \bigg(\int_{0}^{t} u \bigg)^{{1} / {r}} \bigg(\int_{t}^{\infty} \bigg(\int_{t}^{x} w \bigg)^{q'} w(x) \bigg(\esup_{y\in [x,\infty)} v(y)^{-1} \bigg)^{q'} \, dx \bigg)^{1 / q'}.	
    \end{align*}

    Moreover, if $C$ is the best constant in \eqref{main}, then $C \approx F_1 + F_4$.	
    
    {\rm (d)} Let $\max\{q,\,r\} < 1$. Then inequality \eqref{main} holds if and only if $F_3 < \infty$ and $F_5 < \infty$, where
    \begin{align*}
    F_5 & := \bigg(\int_{0}^{\infty}\bigg(\int_{0}^{t} u \bigg)^{r'} u(t) \bigg( \int_{t}^{\infty} \bigg(\int_{t}^{x} w \bigg)^{q'} w(x) \bigg( \esup_{y\in [x,\infty)} v(y)^{-1} \bigg)^{q'} \, dx \bigg)^{r' / q'} \, dt \bigg)^{\frac{1}{r'}}.	
    \end{align*}
    
    Moreover, if $C$ is the best constant in \eqref{main}, then $C \approx F_3 + F_5$.	
\end{thm}

Note that Theorem \ref{main1} coincides with \cite[Theorem 1.3]{krepick}, which was presented with more general weight functions without proof. 

The corrigendum is organized as follows. We start with formulations of background material in Section \ref{pre}. Discrete characterization of inequality {\color{red}{(4.4)}} from \cite{mus.2017} is presented in Section \ref{s.5} (see Lemma \ref{lem.1.6}). Continuous sufficient conditions for inequality {\color{red}{(4.4)}} from \cite{mus.2017} are studied in Section \ref{s.6} (see Lemma \ref{lem.1.8}). In Section \ref{s.7}, we show that these conditions are necessary for inequality \eqref{main} (see Lemma \ref{lem.1.9}). The proof of the main theorem is given in Section \ref{s.8}.


\section{Background Materials}\label{pre}


We recall the definition and some properties of the discretization sequence from \cite{mus.2017}.
\begin{defi}\label{rem.disc.}
	Assume that $u$ is a weight function on $(0,+\infty)$. If $\int_0^{\infty} u(t)\,dt = +\infty$, let $\{x_k\}_{-\infty}^{+\infty} \subset (0,\infty)$ be a strictly increasing sequence  such that $\int_0^{x_k} u(t)\,dt = 2^k$, $k \in \Z$. Denote $M:= +\infty$ and $\zstroke = \Z$, when $\int_0^{\infty} u(t)\,dt = \infty$.  If $\int_0^{\infty} u(t)\,dt < +\infty$, define a strictly increasing sequence $\{x_k\}_{k = -\infty}^M$ such that $\int_0^{x_k} u(t)\,dt = 2^k$, $-\infty < k \le M$, where $M$ satisfies the inequality $2^M \le \int_0^{\infty} u(t)\,dt < 2^{M+1}$. Denote $x_{M+1} : = \infty$ and $\zstroke : = \{k \in \Z:\, k \le M\}$, when $\int_0^{\infty} u(t)\,dt < \infty$.
	Obviously, $\bigcup_{k \in \zstroke} [x_k,x_{k+1}) = (0,\infty)$ in both cases. The sequence $\{x_k\}_{k \in \zstroke}$ is called a covering sequence. 
\end{defi}

\begin{rem}
	We shall use the following equivalences without mentioning anytime we need them.
	
	Assume that $\{x_k\}_{k \in \zstroke}$ is a covering sequence. Clearly, 
	$$
	\int_{x_{k-1}}^{x_k} u \approx 2^k, \quad k \in \zstroke.
	$$
	Moreover,
	$$
	\int_{x_{k-1}}^{x_k} \bigg( \int_{x_{k-1}}^t u \bigg)^{r'} u(t)\,dt \approx \int_{x_{k-1}}^{x_k} \bigg( \int_0^t u \bigg)^{r'} u(t)\,dt  \approx 2^{k(r' + 1)}
	$$
	when $0 < r < 1$.
\end{rem}

The following reduction statement is true.
\begin{thm}\label{thm.1.5}
	Let $0 < q,\, r < \infty$ and  $u,\,v,\,w \in {\mathcal W}\I$.	Assume that $\{x_k\}_{k \in \zstroke}$ is a covering sequence.
	Then inequality \eqref{main} holds if and only if both inequalities 
	\begin{align}
	\bigg\| \bigg\{ 2^{k / r} \bigg(
	\int_{x_k}^{x_{k+1}} \bigg( \int_s^{x_{k+1}} h \bigg)^q
	w(s)\,ds\bigg)^{ 1 / q} \bigg\} \bigg\|_{\ell^r(\zstroke)} & \le C \, \bigg\| \bigg\{ \bigg( \int_{x_n}^{x_{n+1}} h v \bigg) \bigg\} \bigg\|_{\ell^1 (\zstroke)}, \quad h \in {\mathfrak M}^+(0,\infty) \label{ineq.discr} \\
	\intertext{and}
	\bigg( \int_0^{\infty}  u(t) \bigg[ \sup_{s \in [t,\infty)} \bigg(
	\int_t^s w \bigg)^{1 / q} \bigg( \int_s^{\infty} h\bigg) \bigg]^r \,dt\bigg)^{1/r}  & \leq C \, \int_0^{\infty} h v, \quad h \in {\mathfrak M}^+(0,\infty) \label{ineq.supr}
	\end{align}
	hold.
\end{thm}
\begin{proof}
The statement immediately follows from \cite[Lemma 4.5]{mus.2017}. 
\end{proof}	

We recall the following statement concerning the characterization of inequality \eqref{ineq.supr}.
\begin{thm}\cite[Corollary 4.7]{mus.2017}\label{cor.1.1}
	Let $0 < r < \infty$, $0 < q < \infty$, and  $u,\,v,\,w \in {\mathcal W}\I$.
	
	{\rm (a)} Let $r \ge 1$. Then inequality \eqref{ineq.supr} holds if and only if
	\begin{align*}
	D_1 : & = \esup_{t \in (0,\infty)} \bigg( \int_0^t   u(x) \bigg( \int_x^t w \bigg)^{r / q}  \,dx \bigg)^{1 / r} \esup_{s \in [t,\infty)} v(s)^{-1} < \infty, \\
	D_2 : & = \esup_{t \in (0,\infty)}  \bigg( \int_0^t u \bigg)^{1/r}  \bigg( \esup_{s \in [t,\infty)} \bigg( \int_t^s w \bigg)^{1/q} v(s)^{-1} \bigg) < \infty.
	\end{align*}
	Moreover, if $C$ is the best constant in \eqref{ineq.supr}, then $C \approx D_1 + D_2$.
	
	{\rm (b)} Let $r < 1$. Then inequality \eqref{ineq.supr} holds if and only if
	\begin{align*}
	E_1 : & = \bigg( \int_0^{\infty}  \bigg( \int_0^t u \bigg)^{r'}  u(t) \bigg( \esup_{s \in [t,\infty)} \bigg(
	\int_t^s w \bigg)^{r' / q} v(s)^{-r'}\bigg) \,dt\bigg)^{1/r'} < \infty, \\
	E_2 : & = \bigg( \int_0^{\infty}  \bigg( \int_0^t u(x) \bigg( \int_x^t w \bigg)^{r / q} dx \bigg)^{r'} u(t) \bigg( \esup_{s \in [t,\infty)} \bigg( \int_t^s w \bigg)^{r/q} v(s)^{-r'} \bigg) \,dt\bigg)^{1/r'}  < \infty.
	\end{align*}
	Moreover, if $C$ is the best constant in \eqref{ineq.supr}, then $C \approx E_1 + E_2$.	
\end{thm}

\begin{rem}\label{cor.2.0}
	Recall that, if $F$ is a non-negative non-decreasing function on $\I$, then
	\begin{equation}\label{Fubini.2}	
	\esup_{t \in (0,\infty)} F(t)G(t) = \esup_{t \in (0,\infty)} F(t)
	\esup_{\tau \in [t,\infty)} G(\tau)
	\end{equation}
	(see, for instance, \cite{gp2}).
	
	In view of \eqref{Fubini.2}, we have that
	$$
	\esup_{s \in [t,\infty)} \bigg( \int_t^s w \bigg)^{1/q} v(s)^{-1} = \esup_{s \in [t,\infty)} \bigg( \int_t^s w \bigg)^{1/q} \esup_{y \in [s,\infty)} v(y)^{-1}.
	$$
\end{rem}

\begin{rem}\label{cor.2.1}
	Note that $D_2 \lesssim D_1$. Indeed: By Remark \ref{cor.2.0}, we get that
	\begin{align*}
	D_2 = & \esup_{t \in (0,\infty)} \bigg( \int_0^t u \bigg)^{1 / r} \bigg( \esup_{s \in [t,\infty)} \bigg( \int_t^s w \bigg)^{{1} / {q}} \bigg( \esup_{y \in [s,\infty)} v(y)^{-1}\bigg) \bigg) \\
	= & \esup_{t \in (0,\infty)} \bigg( \esup_{s \in [t,\infty)} \bigg( \int_0^t u \bigg)^{1 / r} \bigg( \int_t^s w \bigg)^{{1} / {q}} \bigg( \esup_{y \in [s,\infty)} v(y)^{-1}\bigg) \bigg) \\
	\le & \esup_{t \in (0,\infty)} \bigg( \esup_{s \in [t,\infty)} \bigg( \int_0^t u(x) \bigg( \int_x^s w \bigg)^{r / {q}} \,dx \bigg)^{1 / r}  \bigg( \esup_{y \in [s,\infty)} v(y)^{-1}\bigg) \bigg) \\
	\le & \esup_{t \in (0,\infty)} \bigg( \esup_{s \in [t,\infty)} \bigg( \int_0^s u(x) \bigg( \int_x^s w \bigg)^{r / {q}} \,dx \bigg)^{1 / r}  \bigg( \esup_{y \in [s,\infty)} v(y)^{-1}\bigg) \bigg) \\
	= & \esup_{t \in (0,\infty)} \bigg( \int_0^t u(x) \bigg( \int_x^t w \bigg)^{r / {q}} \,dx \bigg)^{1 / r}  \bigg( \esup_{y \in [t,\infty)} v(y)^{-1}\bigg) \bigg) = D_1.
	\end{align*}
\end{rem}

In view of Remarks \ref{cor.2.0} and \ref{cor.2.1}, Theorem \ref{cor.1.1} reads as follows:
\begin{thm}\label{cor.1.2}
	Let $0 < r < \infty$, $0 < q < \infty$, and  $u,\,v,\,w \in {\mathcal W}\I$.
	
	{\rm (a)} Let $r \ge 1$. Then inequality \eqref{ineq.supr} holds if and only if
	\begin{align*}
	D : & = \esup_{t \in (0,\infty)} \bigg( \int_0^t   u(x) \bigg( \int_x^t w \bigg)^{r / q}  \,dx \bigg)^{1 / r} \bigg( \esup_{y \in [t,\infty)} v(y)^{-1} \bigg) < \infty.
	\end{align*}
	Moreover, if $C$ is the best constant in \eqref{ineq.supr}, then $C \approx D$.
	
	{\rm (b)} Let $r < 1$. Then inequality \eqref{ineq.supr} holds if and only if
	\begin{align*}
	E_1 : & = \bigg( \int_0^{\infty}  \bigg( \int_0^t u \bigg)^{r'}  u(t) \bigg( \esup_{s \in [t,\infty)} \bigg(
	\int_t^s w \bigg)^{r' / q} \bigg( \esup_{y \in [s,\infty)} v(y)^{-1} \bigg)^{r'}\bigg) \,dt\bigg)^{1/r'} < \infty, \\
	E_2 : & = \bigg( \int_0^{\infty}  \bigg( \int_0^t u(x) \bigg( \int_x^t w \bigg)^{r / q} dx \bigg)^{r'} u(t) \bigg( \esup_{s \in [t,\infty)} \bigg( \int_t^s w \bigg)^{r/q} \bigg( \esup_{y \in [s,\infty)} v(y)^{-1} \bigg)^{r'} \bigg) \,dt\bigg)^{1/r'}  < \infty.
	\end{align*}
	Moreover, if $C$ is the best constant in \eqref{ineq.supr}, then $C \approx E_1 + E_2$.	
\end{thm}


\section{Discrete solution of inequality \eqref{ineq.discr}.}\label{s.5}

In this section, we will use the following well-known characterizations of weights for which the weighted Hardy-type inequality holds (see, for instance, \cite{ok}).		
\begin{thm}\label{thm.Copson}
	Let $0 < q < \infty$ and $v,\,w \in {\mathcal W} (0,\infty)$.
	
	{\rm (i)} Let $1 \le q$. Then 
	$$
	\sup_{f \in \mp^+ (0,\infty)} \frac{\bigg( \int_0^{\infty} \bigg( \int_x^{\infty} f \bigg)^q w(x)\,dx\bigg)^{{1} / {q}}}{\int_0^{\infty} f v} \approx \sup_{x \in (0,\infty)} \bigg( \int_0^x w \bigg)^{{1} / {q}} \bigg( \esup_{t \in [x,\infty)} v(t)^{-1}\bigg).
	$$
	
	{\rm (ii)} Let $q < 1$. Then 
	$$
	\sup_{f \in \mp^+ (0,\infty)} \frac{\bigg( \int_0^{\infty} \bigg( \int_x^{\infty} f \bigg)^q w(x)\,dx\bigg)^{{1} / {q}}}{\int_0^{\infty} f v} \approx \bigg( \int_0^{\infty} \bigg( \int_0^x w \bigg)^{q'} \, w(x) \, \bigg( \esup_{t \in [x,\infty)} v(t)^{-1}\bigg)^{q'} \,dx \bigg)^{{1} / {q'}}.
	$$
\end{thm}	

Now we give a discrete characterization of inequality \eqref{ineq.discr}. Recall that $\rho$ is defined for any $0 < r < \infty$ in \cite{mus.2017} by $1 / \rho : = ( 1 / r - 1)_+$, where $a_+ = a$ when $a>0$ and $a_+ = 0$ when $a \le 0$.
\begin{lem}\label{lem.1.6}
	Let $0 < q,\, r < \infty$ and  $u,\,w,\,v \in {\mathcal W}(0,\infty)$. Assume that $\{x_k\}_{k \in \zstroke}$ is a covering sequence. 
	
	{\rm (i)} If $q \ge 1$, then inequality \eqref{ineq.discr}
	holds with constant independent of $h \in {\mathfrak M}^+(0,\infty)$ if and only if $A_1 < \infty$, where
	$$
	A_1 := \bigg\| \bigg\{ 2^{k / r} \bigg( \esup_{x \in [x_k,x_{k+1})} \bigg( \int_{x_k}^x w \bigg)^{{1} / {q}} \bigg( \esup_{t \in [x,x_{k+1})} v(t)^{-1}\bigg) \bigg) \bigg\} \bigg\|_{\ell^{\rho}(\zstroke)}.
	$$
	Moreover, if $C$ is the best constant in \eqref{ineq.discr}, then $C \approx A_1$.	
	
	{\rm (ii)} If $q < 1$, then inequality \eqref{ineq.discr}
	holds with constant independent of $h \in {\mathfrak M}^+(0,\infty)$ if and only if $A_2 < \infty$, where
	$$
	A_2 := \bigg\| \bigg\{ 2^{k / r} \bigg( \int_{x_k}^{x_{k+1}} \bigg( \int_{x_k}^x w \bigg)^{q'} w(x) \bigg( \esup_{t \in [x,x_{k+1})} v(t)^{-1} \bigg)^{q'} \,dx\bigg)^{ 1 / {q'}} \bigg\} \bigg\|_{\ell^{\rho}(\zstroke)}.
	$$
	Moreover, if $C$ is the best constant in \eqref{ineq.discr}, then $C \approx A_2$.	
\end{lem}

\begin{proof}
We give the proof of the second case. The first one can be proved similiarily, so its proof is omitted. 
	
	 {\bf Necessity.} Let $q < 1$. Assume that inequality \eqref{ineq.discr}
	 holds with constant $C$ independent of $h \in {\mathfrak M}^+(0,\infty)$. 
	 
	 Since, for any $k \in \zstroke$
	 $$
	 \sup \bigg\{ \int_{x_k}^{x_{k+1}} h(s)\,ds :\,  \int_{x_k}^{x_{k+1}} hv \le 1 \bigg\} = \esup_{t \in [x_k,x_{k+1})} v(t)^{-1},
	 $$
	 then there exists $h_k \in \mp^+(0,\infty)$ such that $\supp h_k \in (x_k,x_{k+1})$, $\int_{x_k}^{x_{k+1}} h_k v = 1$ and
	 $$
	 \int_{x_k}^{x_{k+1}} h(s)\,ds \ge \frac{1}{2} \esup_{t \in [x_k,x_{k+1})} v(t)^{-1}.
	 $$
	 
	 Define
	 \begin{equation*}
	 h = \sum_{m \in \zstroke} a_m h_m,
	 \end{equation*}
	 where $\{a_k\}_{k \in \zstroke}$ is any sequence of nonnegative numbers such that $\sum_{k \in \zstroke} a_k < \infty$.
	 
	 Obviously, $h \in \mp^+(0,\infty)$, and inequality \eqref{ineq.discr} holds with $h$. Consequently, the inequality
	 \begin{equation*}
	 \bigg\| \bigg\{ a_k 2^{k / {r}} \bigg( \int_{x_k}^{x_{k+1}} \bigg( \int_{x_k}^x w \bigg)^{q'} w(x) \bigg( \esup_{t \in [x,x_{k+1})} v(t)^{-1} \bigg)^{q'} \,dx \bigg)^{{1} / {q'}} \bigg\} \bigg\|_{\ell^r (\zstroke)} \le C \, \|\{a_k\}\|_{\ell^1(\zstroke)}
	 \end{equation*}
	 holds for all sequences $\{a_k\}_{k \in \zstroke}$ of nonnegative numbers. 
	 
	 By \cite[Proposition 2.4]{mus.2017}, we arrive at
	 $$
	 A_2 = \bigg\| \bigg\{ 2^{k / {r}} \bigg( \int_{x_k}^{x_{k+1}} \bigg( \int_{x_k}^x w \bigg)^{q'} w(x) \bigg( \esup_{t \in [x,x_{k+1})} v(t)^{-1} \bigg)^{q'} \,dx \bigg)^{{1} / {q'}} \bigg\} \bigg\|_{\ell^{\rho} (\zstroke)}
	 \le C.
	 $$
	 
	{\bf Sufficiency.}	Assume that $A_2 < \infty$.	By Theorem \ref{thm.Copson}, (ii), applying the discrete H\"{o}lder inequality, we have that
	\begin{align*}
	\bigg\| \bigg\{ 2^{k / r} \bigg( \int_{x_k}^{x_{k+1}} \bigg( \int_s^{x_{k+1}} h \bigg)^q
	w(s)\,ds\bigg)^{ 1 / q} \bigg\} \bigg\|_{\ell^r(\zstroke)} & \\
	& \hspace{-3cm} \le \bigg\| \bigg\{ 2^{k / r} \bigg( \int_{x_k}^{x_{k+1}} \bigg( \int_{x_k}^x w \bigg)^{q'} w(x) \bigg( \esup_{t \in [x,x_{k+1})} v(t)^{-1} \bigg)^{q'} \,dx\bigg)^{ 1 / {q'}} \bigg( \int_{x_k}^{x_{k+1}} h v\bigg)\bigg\} \bigg\|_{\ell^r(\zstroke)} \\
	& \hspace{-3cm} \le 
	\bigg\| \bigg\{ 2^{k / r} \bigg( \int_{x_k}^{x_{k+1}} \bigg( \int_{x_k}^x w \bigg)^{q'} w(x) \bigg( \esup_{t \in [x,x_{k+1})} v(t)^{-1} \bigg)^{q'} \,dx\bigg)^{ 1 / {q'}} \bigg\|_{\ell^{\rho}(\zstroke)}  \bigg\| \bigg\{ \int_{x_n}^{x_{n+1}} h v \bigg\} \bigg\|_{\ell^1 (\zstroke)} \\
	& \hspace{-3cm} = A_2 \, \bigg\| \bigg\{ \int_{x_n}^{x_{n+1}} h v \bigg\} \bigg\|_{\ell^1 (\zstroke)}.
	\end{align*} 
	
	Thus, inequality \eqref{ineq.discr} holds, and, if $C$ is the best constant in \eqref{ineq.discr}, then 		
	$$
	C \le A_2.
	$$	
	
	The proof is completed.
\end{proof}


\section{Continuous sufficient conditions for inequality \eqref{ineq.discr}}\label{s.6}

\begin{lem}\label{lem.1.8}
	Let $0 < q,\, r < \infty$ and  $u,\,v,\,w \in {\mathcal W}(0,\infty)$. 
	
	{\rm (i)} Let $1 \le \min\{q,\,r\}$. If
	$$
	B_1 : = \esup_{t \in (0,\infty)} \bigg( \int_0^t u \bigg)^{1 / r} \bigg( \esup_{x \in [t,\infty)} \bigg( \int_t^x w \bigg)^{{1} / {q}} \bigg( \esup_{y \in [x,\infty)} v(y)^{-1}\bigg) \bigg) < \infty,
	$$
	then inequality \eqref{ineq.discr}	holds with constant independent of $h \in {\mathfrak M}^+(0,\infty)$. Moreover, if $C$ is the best constant in \eqref{ineq.discr}, then $C \lesssim B_1$.	

    {\rm (ii)} Let $r < 1 \le q$. If
    $$
    B_2 : =  \bigg( \int_0^{\infty} \bigg( \int_0^t u \bigg)^{r'} u(t) \bigg( \esup_{x \in [t,\infty)} \bigg( \int_t^x w \bigg)^{{1} / {q}} \bigg( \esup_{y \in [x,\infty)} v(y)^{-1}\bigg) \bigg)^{r'} \,dt \bigg)^{1 / {r'}} < \infty,
    $$
    then inequality \eqref{ineq.discr}	holds with constant independent of $h \in {\mathfrak M}^+(0,\infty)$. Moreover, if $C$ is the best constant in \eqref{ineq.discr}, then $C \lesssim B_2$.	
    
	{\rm (iii)} Let $q < 1 \le r$. If
    $$
    B_3 : = \sup_{t \in (0,\infty)} \bigg( \int_0^t u \bigg)^{1 / r} \bigg( \int_t^{\infty} \bigg( \int_t^x w \bigg)^{q'} w(x) \bigg( \esup_{y \in [x,\infty)} v(y)^{-1} \bigg)^{q'} \,dx\bigg)^{ 1 / {q'}} < \infty,
    $$
    then inequality \eqref{ineq.discr}	holds with constant independent of $h \in {\mathfrak M}^+(0,\infty)$. Moreover, if $C$ is the best constant in \eqref{ineq.discr}, then $C \lesssim B_3$.	
	
	{\rm (iv)} Let $\max\{q,\,r\} < 1$. If
	$$
	B_4 : = \bigg( \int_0^{\infty} \bigg( \int_0^t u \bigg)^{r'} u(t) \bigg( \int_t^{\infty} \bigg( \int_t^x w \bigg)^{q'} w(x) \bigg( \esup_{y \in [x,\infty)} v(y)^{-1} \bigg)^{q'} \,dx\bigg)^{ {r'} / {q'}}\,dt \bigg)^{ 1 / {r'}} < \infty,
	$$
	then inequality \eqref{ineq.discr}	holds with constant independent of $h \in {\mathfrak M}^+(0,\infty)$. Moreover, if $C$ is the best constant in \eqref{ineq.discr}, then $C \lesssim B_4$.	
\end{lem}
	
\begin{proof}
{\rm (i)} Let $1 \le \min\{q,\,r\}$. Assume that $B_1 < \infty$. Clearly,
\begin{align*}
A_1 & = \sup_{k \in \zstroke} 2^{k / r} \bigg( \esup_{x \in [x_k,x_{k+1})} \bigg( \int_{x_k}^x w \bigg)^{{1} / {q}} \bigg( \esup_{y \in [x,x_{k+1})} v(y)^{-1}\bigg) \bigg) \\
& \lesssim \sup_{k \in \zstroke} \bigg( \int_0^{x_k} u \bigg)^{1 / r} \bigg( \esup_{x \in [x_k,\infty)} \bigg( \int_{x_k}^x w \bigg)^{{1} / {q}} \bigg( \esup_{y \in [x,\infty)} v(y)^{-1}\bigg) \bigg) \\
& \le \sup_{t \in (0,\infty)} \bigg( \int_0^t u \bigg)^{1 / r} \bigg( \esup_{x \in [t,\infty)} \bigg( \int_t^x w \bigg)^{{1} / {q}} \bigg( \esup_{y \in [x,\infty)} v(y)^{-1}\bigg) \bigg) = B_1.
\end{align*}
Thus $A_1 < \infty$, and the statement follows from Lemma \ref{lem.1.6}, (i).

{\rm (ii)} Let $r < 1 \le q$. Assume that $B_2 < \infty$. Clearly,
\begin{align*}
A_1 & = \bigg( \sum_{k \in \zstroke} 2^{k (r'+ 1)} \bigg( \esup_{x \in [x_k,x_{k+1})} \bigg( \int_{x_k}^x w \bigg)^{{1} / {q}} \bigg( \esup_{y \in [x,x_{k+1})} v(y)^{-1}\bigg) \bigg)^{r'} \bigg)^{1 / {r'}} \\
& \lesssim \bigg( \sum_{k \in \zstroke} \int_{x_{k-1}}^{x_k} \bigg( \int_0^t u \bigg)^{r'} u(t)\,dt  \cdot \bigg( \esup_{x \in [x_k,\infty)} \bigg( \int_{x_k}^x w(s)\,ds\bigg)^{{1} / {q}} \bigg( \esup_{y \in [x,\infty)} v(y)^{-1}\bigg) \bigg)^{r'} \bigg)^{1 / {r'}} \\
& \le \bigg( \sum_{k \in \zstroke} \int_{x_{k-1}}^{x_k} \bigg( \int_0^t u \bigg)^{r'} u(t) \bigg( \esup_{x \in [t,\infty)} \bigg( \int_t^x w \bigg)^{{1} / {q}} \bigg( \esup_{y \in [x,\infty)} v(y)^{-1}\bigg) \bigg)^{r'} \,dt \bigg)^{1 / {r'}} \\
& \le \bigg( \int_0^{\infty} \bigg( \int_0^t u \bigg)^{r'} u(t) \bigg( \esup_{x \in [t,\infty)} \bigg( \int_t^x w \bigg)^{{1} / {q}} \bigg( \esup_{y \in [x,\infty)} v(y)^{-1}\bigg) \bigg)^{r'} \,dt \bigg)^{1 / {r'}} = B_2.
\end{align*}
Thus $A_1 < \infty$, and the statement follows from Lemma \ref{lem.1.6}, (i).

{\rm (iii)} Let $q < 1 \le r$. Assume that $B_3 < \infty$. We have that
\begin{align*}
A_2 & = \sup_{k \in \zstroke} 2^{k / r} \bigg( \int_{x_k}^{x_{k+1}} \bigg( \int_{x_k}^x w \bigg)^{q'} w(x) \bigg( \esup_{y \in [x,x_{k+1})} v(y)^{-1} \bigg)^{q'} \,dx\bigg)^{ 1 / {q'}}  \\
& \lesssim \sup_{k \in \zstroke} \bigg( \int_0^{x_k} u \bigg)^{1 / r} \bigg( \int_{x_k}^{\infty} \bigg( \int_{x_k}^x w \bigg)^{q'} w(x) \bigg( \esup_{y \in [x,\infty)} v(y)^{-1} \bigg)^{q'} \,dx\bigg)^{ 1 / {q'}}  \\
& \le \sup_{t \in (0,\infty)} \bigg( \int_0^t u \bigg)^{1 / r} \bigg( \int_t^{\infty} \bigg( \int_t^x w \bigg)^{q'} w(x) \bigg( \esup_{y \in [x,\infty)} v(y)^{-1} \bigg)^{q'} \,dx\bigg)^{ 1 / {q'}} = B_3.
\end{align*}
Thus $A_2 < \infty$, and the statement follows from Lemma \ref{lem.1.6}, (ii).	

{\rm (iv)} Let $\max\{q,\,r\} < 1$. Assume that $B_4 < \infty$. We have that
\begin{align*}
A_2 & = \bigg( \sum_{k \in \zstroke} 2^{k (r' + 1)} \bigg( \int_{x_k}^{x_{k+1}} \bigg( \int_{x_k}^x w \bigg)^{q'} w(x) \bigg( \esup_{y \in [x,x_{k+1})} v(y)^{-1} \bigg)^{q'} \,dx\bigg)^{ {r'} / {q'}} \bigg)^{ 1 / {r'}} \\
& \lesssim \bigg( \sum_{k \in \zstroke} \int_{x_{k-1}}^{x_k} \bigg( \int_0^t u \bigg)^{r'} u(t)\,dt  \cdot \bigg( \int_{x_k}^{\infty} \bigg( \int_{x_k}^x w \bigg)^{q'} w(x) \bigg( \esup_{y \in [x,\infty)} v(y)^{-1} \bigg)^{q'} \,dx \bigg)^{ {r'} / {q'}} \bigg)^{ 1 / {r'}} \\
& \le \bigg( \sum_{k \in \zstroke} \int_{x_{k-1}}^{x_k} \bigg( \int_0^t u \bigg)^{r'} u(t) \bigg( \int_t^{\infty} \bigg( \int_t^x w \bigg)^{q'} w(x) \bigg( \esup_{y \in [x,\infty)} v(y)^{-1} \bigg)^{q'} \,dx\bigg)^{ {r'} / {q'}}\,dt \bigg)^{ 1 / {r'}} \\
& \le \bigg( \int_0^{\infty} \bigg( \int_0^t u \bigg)^{r'} u(t) \bigg( \int_t^{\infty} \bigg( \int_t^x w \bigg)^{q'} w(x) \bigg( \esup_{y \in [x,\infty)} v(y)^{-1} \bigg)^{q'} \,dx\bigg)^{ {r'} / {q'}}\,dt \bigg)^{ 1 / {r'}} = B_4. 
\end{align*}
Thus $A_2 < \infty$, and the statement follows from Lemma \ref{lem.1.6}, (ii).

The proof is completed.
\end{proof}	


\section{Necessity of conditions $B_i$, $i=\overline{1,4}$ for inequality \eqref{main}}\label{s.7}

We need the following statement.
\begin{lem}\label{lem.2.3}
	Let $0 < q,\, r < \infty$ and  $u,\,v,\,w \in {\mathcal W}\I$. Assume that inequality \eqref{main}. Then
	\begin{align*}
	S : = \sup_{\{b_k\}} \bigg\| \bigg\{ \bigg( \int_{b_{k-1}}^{b_k} \bigg( \int_x^{b_k} w \bigg)^{{r} / {q}} u(x)\,dx \bigg)^{1 / r}  \bigg( \esup_{t \in [b_k,b_{k+1})} v(t)^{-1} \bigg) \bigg\} \bigg\|_{\ell^{\rho} (\Z)} < \infty.
	\end{align*}
	Moreover, if $C$ is the best constant in \eqref{main}, then $S \lesssim C$.
	
	Here the suprema are taken over all increasing sequences $\{b_k\} = \{b_k\}_{k \in \Z} \subset (0,\infty)$.
\end{lem}

\begin{proof}
	Assume that inequality \eqref{main} holds. We take any increasing sequence $\{b_k\} = \{b_k\}_{k \in \Z} \subset (0,\infty)$
	
	Since, for any $k \in \Z$
	$$
	\sup \bigg\{ \int_{b_k}^{b_{k+1}} h(s)\,ds :\,  \int_{b_k}^{b_{k+1}} hv \le 1 \bigg\} = \esup_{t \in [b_k,b_{k+1})} v(t)^{-1},
	$$
	then there exists $h_k \in \mp^+(0,\infty)$ such that $\supp h_k \in (b_k,b_{k+1})$, $\int_{b_k}^{b_{k+1}} h_k v = 1$ and
	$$
	\int_{b_k}^{b_{k+1}} h(s)\,ds \ge \frac{1}{2} \esup_{t \in [b_k,b_{k+1})} v(t)^{-1}.
	$$
	
	Define
	\begin{equation*}
	h = \sum_{m \in \Z} a_m h_m,
	\end{equation*}
	where $\{a_k\}_{k \in \Z}$ is any sequence of nonnegative numbers such that $\sum_{k \in \Z} a_k < \infty$. Clearly, $h \in \mp^+(0,\infty)$, and inequality \eqref{main} holds with $h$.
	
	Writing
	\begin{align*}
	\bigg( \int_0^{\infty} \bigg( \int_x^{\infty} \bigg( \int_t^{\infty} h \bigg)^q w(t)\,dt
	\bigg)^{r / q} u(x)\,dx \bigg)^{1/r} & \\
	& \hspace{-5cm} \ge  \bigg( \sum_{k = - \infty}^{\infty} \int_{b_{k-1}}^{b_k} \bigg( \int_x^{b_k} \bigg( \int_t^{b_{k + 1}} h \bigg)^q w(t)\,dt \bigg)^{{r} / {q}} u(x)\,dx \bigg)^{\frac{1}{r}} \\
	& \hspace{-5cm} \ge  \bigg( \sum_{k = - \infty}^{\infty} \int_{b_{k-1}}^{b_k} \bigg( \int_x^{b_k} w \bigg)^{{r} / {q}} u(x)\,dx \cdot \bigg( \int_{b_k}^{b_{k + 1}} h \bigg)^r \bigg)^{\frac{1}{r}} \\
	& \hspace{-5cm} =  \bigg( \sum_{k = - \infty}^{\infty} \bigg( \bigg( \int_{b_{k-1}}^{b_k} \bigg( \int_x^{b_k} w \bigg)^{{r} / {q}} u(x)\,dx \bigg)^{1 / r}  \bigg( \int_{b_k}^{b_{k + 1}} h \bigg) \bigg)^r \bigg)^{\frac{1}{r}} \\
	& \hspace{-5cm} =  \bigg( \sum_{k = - \infty}^{\infty} \bigg( a_k \bigg( \int_{b_{k-1}}^{b_k} \bigg( \int_x^{b_k} w \bigg)^{{r} / {q}} u(x)\,dx \bigg)^{1 / r}  \bigg( \int_{b_k}^{b_{k + 1}} h_k \bigg) \bigg)^r \bigg)^{\frac{1}{r}} \\
	& \hspace{-5cm} \gtrsim  \bigg( \sum_{k = - \infty}^{\infty} \bigg( a_k \bigg( \int_{b_{k-1}}^{b_k} \bigg( \int_x^{b_k} w \bigg)^{{r} / {q}} u(x)\,dx \bigg)^{1 / r}  \bigg( \esup_{t \in [b_k,b_{k+1})} v(t)^{-1} \bigg) \bigg)^r \bigg)^{\frac{1}{r}}
	\end{align*}
	and
	\begin{align*}
	\int_0^{\infty} hv = \sum_{k = - \infty}^{\infty} \int_{b_k}^{b_{k+1}} hv = \sum_{k = - \infty}^{\infty} a_k \int_{b_k}^{b_{k+1}} h_k v = \sum_{k = - \infty}^{\infty} a_k,
	\end{align*}
	we get that the inequality
	\begin{equation*}\label{disc.ineq.2.2}
	\bigg\| \bigg\{ a_k \bigg( \int_{b_{k-1}}^{b_k} \bigg( \int_x^{b_k} w \bigg)^{{r} / {q}} u(x)\,dx \bigg)^{1 / r}  \bigg( \esup_{t \in [b_k,b_{k+1})} v(t)^{-1} \bigg) \bigg\} \bigg\|_{\ell^r (\Z)} \le C \, \|\{a_k\}\|_{\ell^1(\Z)}
	\end{equation*}
	holds for all sequences $\{a_k\}_{k \in \Z}$ of nonnegative numbers. 
	
	By \cite[Proposition 2.4]{mus.2017}, we arrive at
	$$
	\bigg\| \bigg\{ \bigg( \int_{b_{k-1}}^{b_k} \bigg( \int_x^{b_k} w \bigg)^{{r} / {q}} u(x)\,dx \bigg)^{1 / r}  \bigg( \esup_{t \in [b_k,b_{k+1})} v(t)^{-1} \bigg) \bigg\} \bigg\|_{\ell^{\rho} (\Z)} \le C.
	$$
	
	Consequently,
	$$
	S \le C.
	$$

    The proof is completed.
\end{proof}	

In order to prove the following statement we use the blocking technique developed in \cite{grosse} and \cite{GogStep}.
\begin{lem}\label{lem.2.5}
	Let $0 < q < 1$, $0 < r < \infty$ and  $u,\,v,\,w \in {\mathcal W}\I$. Assume that $\{x_k\}_{k \in \zstroke}$ is a covering sequence. Suppose that $S < \infty$. Then
	\begin{align*}
	\bigg\| \bigg\{ 2^{{k} / {r}} \bigg( \int_{x_k}^{x_{k + 1}} w \bigg) \bigg( \int_{x_k}^{\infty} w(x) \bigg( \esup_{t \in [x,\infty)} v(t)^{-1}\bigg)^{q'}\,dx \bigg)^{1 / {q'}} \bigg\} \bigg\|_{\ell^{\rho}(\zstroke \backslash \{M\})} \le S.
	\end{align*}
\end{lem}

\begin{proof}
{\rm (i)} Let $r < 1$. For any $i \in \Z$ we set 
\begin{equation*}
{\mathcal N}_i^1 : = \left\{ m \in \zstroke:\, 2^i \le \sum_{k < m} 2^{k(r' + 1)} \bigg( \int_{x_k}^{x_{k+1}} w \bigg)^{r'} < 2^{i + 1} \right\}
\end{equation*}
and denote by ${\mathcal N}^1 : = \bigcup_n {\mathcal N}_{i_n}^1 \subset \Z$, where ${\mathcal N}_{i_n}^1 \neq \emptyset$ and ${\mathcal N}_i^1 = \emptyset$ for $i \not \in \{i_n\}$. Next we define
$$
m_{i_n} = \inf {\mathcal N}_{i_n}^1, \qquad m_{i_n}^+ = \sup {\mathcal N}_{i_n}^1, \qquad z_n = x_{m_{i_n}}.
$$
Assume, for simplicity that $\card \bigcup_n \{ i_n \} = \aleph_0$.

Applying \cite[Lemma 2.3]{mus.2017}, we get that
\begin{align*}
\bigg\| \bigg\{ 2^{{k} / {r}} \bigg( \int_{x_k}^{x_{k + 1}} w \bigg) \bigg( \int_{x_k}^{\infty} w(x) \bigg( \esup_{t \in [x,\infty)} v(t)^{-1}\bigg)^{q'}\,dx \bigg)^{1 / {q'}} \bigg\} \bigg\|_{\ell^{\rho}(\zstroke \backslash \{M\})} & \\
& \hspace{-7cm} = \bigg( \sum_{k < M} 2^{k(r' + 1)} \bigg( \int_{x_k}^{x_{k+1}} w \bigg)^{r'} \bigg( \int_{x_k}^{\infty} w(x) \bigg( \esup_{t \in [x,\infty)} v(t)^{-1} \bigg)^{q'} \,dx \bigg)^{ {r'} / {q'}} \bigg)^{ 1 / {r'}} \\
& \hspace{-7cm} =  \bigg( \sum_{n} \sum_{m_{i_n} \le k \le m_{i_n}^+} 2^{k(r' + 1)} \bigg( \int_{x_k}^{x_{k+1}} w \bigg)^{r'} 
\bigg( \int_{x_k}^{\infty} w(x) \bigg( \esup_{t \in [x,\infty)} v(t)^{-1} \bigg)^{q'} \,dx \bigg)^{ {r'} / {q'}} \bigg)^{ 1 / {r'}} \\
& \hspace{-7cm} \le \bigg( \sum_{n} \bigg( \sum_{m_{i_n} \le k \le m_{i_n}^+} 2^{k(r' + 1)} \bigg( \int_{x_k}^{x_{k+1}} w \bigg)^{r'} \bigg)
\bigg( \int_{z_n}^{\infty} w(x) \bigg( \esup_{t \in [x,\infty)} v(t)^{-1} \bigg)^{q'} \,dx \bigg)^{ {r'} / {q'}} \bigg)^{ 1 / {r'}} \\
& \hspace{-7cm} \lesssim \bigg( \sum_{n} 2^{i_n}
\bigg( \int_{z_n}^{\infty} w(x) \bigg( \esup_{t \in [x,\infty)} v(t)^{-1} \bigg)^{q'} \,dx \bigg)^{ {r'} / {q'}} \bigg)^{ 1 / {r'}} \\
& \hspace{-7cm} = \bigg( \sum_{n} 2^{i_n}
\bigg( \sum_{n \le j} \int_{z_j}^{z_{j+1}} w(x) \bigg( \esup_{t \in [x,\infty)} v(t)^{-1} \bigg)^{q'} \,dx \bigg)^{ {r'} / {q'}} \bigg)^{ 1 / {r'}} \\
& \hspace{-7cm} \approx \bigg( \sum_{n} 2^{i_n}
\bigg( \int_{z_n}^{z_{n+1}} w(x) \bigg( \esup_{t \in [x,\infty)} v(t)^{-1} \bigg)^{q'} \,dx \bigg)^{ {r'} / {q'}} \bigg)^{ 1 / {r'}} \\
& \hspace{-7cm} \le \bigg( \sum_{n} 2^{i_n}
\bigg( \int_{z_n}^{z_{n+1}} w \bigg)^{ {r'} / {q'}} \bigg( \esup_{t \in [z_n,\infty)} v(t)^{-1} \bigg)^{r'} \bigg)^{ 1 / {r'}}.
\end{align*}

For each $j \in \Z$ we set 
\begin{equation*}
{\mathcal N}_j^2 : = \left\{ l \in \Z:\, 2^j \le \sum_{n < l} 2^{i_n} \bigg( \int_{z_n}^{z_{n+1}} w \bigg)^{{r'} / {q'}} < 2^{j + 1} \right\}
\end{equation*}
and denote by ${\mathcal N}^2 : = \bigcup_k {\mathcal N}_{j_k}^2 \subset \Z$, where ${\mathcal N}_{j_k}^2 \neq \emptyset$ and ${\mathcal N}_j^2 = \emptyset$ for $j \not \in \{j_k\}$. Next we define
$$
l_{j_k} = \inf {\mathcal N}_{j_k}^2, \qquad l_{j_k}^+ = \sup {\mathcal N}_{j_k}^2, \qquad y_k = z_{l_{j_k}}.
$$
Suppose, for simplicity that $\card \bigcup_k \{ j_k \} = \aleph_0$.

We have, by \cite[Lemma 2.3]{mus.2017}, that
\begin{align*}
\bigg( \sum_{n} 2^{i_n}
\bigg( \int_{z_n}^{z_{n+1}} w \bigg)^{ {r'} / {q'}} \bigg( \esup_{t \in [z_n,\infty)} v(t)^{-1} \bigg)^{r'} \bigg)^{ 1 / {r'}} & \\
& \hspace{-5cm} = \bigg( \sum_{k} \sum_{l_{j_k} \le n \le l_{j_k}^+} 2^{i_n} \bigg( \int_{z_n}^{z_{n+1}} w \bigg)^{{r'} / {q'}} \bigg( \esup_{t \in [z_n,\infty)} v(t)^{-1} \bigg)^{r'} \bigg)^{ 1 / {r'}} \\
& \hspace{-5cm} \le \bigg( \sum_{k} \bigg( \sum_{l_{j_k} \le n \le l_{j_k}^+} 2^{i_n} \bigg( \int_{z_n}^{z_{n+1}} w \bigg)^{{r'} / {q'}} \bigg) \bigg( \esup_{t \in [y_k,\infty)} v(t)^{-1} \bigg)^{r'} \bigg)^{ 1 / {r'}} \\
& \hspace{-5cm} \lesssim \bigg( \sum_{k} 2^{j_k} \bigg( \esup_{t \in [y_k,\infty)} v(t)^{-1} \bigg)^{r'} \bigg)^{ 1 / {r'}} \\
& \hspace{-5cm} = \bigg( \sum_{k} 2^{j_k} \bigg( \sup_{k \le m} \esup_{t \in [y_m,y_{m+1})} v(t)^{-1} \bigg)^{r'} \bigg)^{ 1 / {r'}} \\
& \hspace{-5cm} \approx \bigg( \sum_{k} 2^{j_k} \bigg( \esup_{t \in [y_k,y_{k+1})} v(t)^{-1} \bigg)^{r'} \bigg)^{ 1 / {r'}}.
\end{align*}

Note that $2^{j_k} \le 2\,\sum_{l_{j_{k-2}} \le n < l_{j_k}} 2^{i_n} \bigg( \int_{z_n}^{z_{n+1}} w \bigg)^{r' / q'}$. Hence
\begin{align*}
\bigg( \sum_{n} 2^{i_n}
\bigg( \int_{z_n}^{z_{n+1}} w \bigg)^{ {r'} / {q'}} \bigg( \esup_{t \in [z_n,\infty)} v(t)^{-1} \bigg)^{r'} \bigg)^{ 1 / {r'}} & \\
& \hspace{-5cm} \lesssim \bigg( \sum_{k} \bigg( \sum_{l_{j_{k-2}} \le n < l_{j_k}} 2^{i_n} \bigg( \int_{z_n}^{z_{n + 1}} w \bigg)^{{r'} / {q'}} \bigg)\bigg( \esup_{t \in [y_k,y_{k+1})} v(t)^{-1} \bigg)^{r'} \bigg)^{ 1 / {r'}}.
\end{align*}

Since $2^{i_n} \le 2\,\sum_{m_{i_{n-2}} \le k < m_{i_n}} 2^{k(r' + 1)} \bigg( \int_{x_k}^{x_{k+1}} w \bigg)^{r'}$, we have that
\begin{align}
2^{i_n} & \lesssim \sum_{m_{i_{n-2}} \le k < m_{i_n}} \bigg( 2^{k} \bigg( \int_{x_k}^{x_{k+1}} w \bigg)^{r} \bigg)^{(r' + 1)} \notag \\
& \le \bigg( \sum_{m_{i_{n-2}} \le k < m_{i_n}} 2^{k} \bigg( \int_{x_k}^{x_{k+1}} w \bigg)^{r} \bigg)^{(r' + 1)} \notag \\
& \approx \bigg( \sum_{m_{i_{n-2}} \le k < m_{i_n}} \int_{x_{k-1}}^{x_k} u(t)\,dt \cdot \bigg( \int_{x_k}^{x_{k+1}} w \bigg)^{r} \bigg)^{(r' + 1)} \notag \\
& \le \bigg( \sum_{m_{i_{n-2}} \le k < m_{i_n}} \int_{x_{k-1}}^{x_k} u(t) \bigg( \int_t^{x_{k+1}} w \bigg)^{r}\,dt \bigg)^{(r' + 1)} \notag \\
& \le \bigg( \int_{z_{n-3}}^{z_n} u(t) \bigg( \int_t^{z_n} w \bigg)^{r}\,dt \bigg)^{(r' + 1)}. \label{eq.1.1.1.1}
\end{align}

On using \eqref{eq.1.1.1.1}, we get that
\begin{align*}
\bigg( \sum_{n} 2^{i_n} \bigg( \int_{z_n}^{z_{n+1}} w \bigg)^{ {r'} / {q'}} \bigg( \esup_{t \in [z_n,\infty)} v(t)^{-1} \bigg)^{r'} \bigg)^{ 1 / {r'}} & \\
& \hspace{-5cm} \lesssim \bigg( \sum_{k} \bigg( \sum_{l_{j_{k-2}} \le n < l_{j_k}} \bigg( \int_{z_{n-3}}^{z_n} u(t) \bigg( \int_t^{z_n} w \bigg)^{r}\,dt \bigg)^{(r' + 1)} \bigg( \int_{z_n}^{z_{n + 1}} w \bigg)^{{r'} / {q'}} \bigg)\bigg( \esup_{t \in [y_k,y_{k+1})} v(t)^{-1} \bigg)^{r'} \bigg)^{ 1 / {r'}} \\
& \hspace{-5cm} \le \bigg( \sum_{k} \bigg( \sum_{l_{j_{k-2}} \le n < l_{j_k}} \bigg( \int_{z_{n-3}}^{z_{n+1}} u(t) \bigg( \int_t^{z_{n+1}} w \bigg)^{r} \bigg( \int_t^{z_{n + 1}} w \bigg)^{r / {q'}}\,dt \bigg)^{(r' + 1)} \bigg)\bigg( \esup_{t \in [y_k,y_{k+1})} v(t)^{-1} \bigg)^{r'} \bigg)^{ 1 / {r'}} \\
& \hspace{-5cm} = \bigg( \sum_{k} \bigg( \sum_{l_{j_{k-2}} \le n < l_{j_k}} \bigg( \int_{z_{n-3}}^{z_{n+1}} u(t) \bigg( \int_t^{z_{n+1}} w \bigg)^{r / q} \,dt \bigg)^{(r' + 1)} \bigg)\bigg( \esup_{t \in [y_k,y_{k+1})} v(t)^{-1} \bigg)^{r'} \bigg)^{ 1 / {r'}} \\
& \hspace{-5cm} \le \bigg( \sum_{k} \bigg( \sum_{l_{j_{k-2}} \le n < l_{j_k}} \int_{z_{n-3}}^{z_{n+1}} u(t) \bigg( \int_t^{z_{n+1}} w \bigg)^{r / q} \,dt \bigg)^{(r' + 1)} \bigg( \esup_{t \in [y_k,y_{k+1})} v(t)^{-1} \bigg)^{r'} \bigg)^{ 1 / {r'}} \\
& \hspace{-5cm} \le \bigg( \sum_{k} \bigg( \int_{y_{k-8}}^{y_k} u(t) \bigg( \int_t^{y_k} w \bigg)^{r / q} \,dt \bigg)^{(r' + 1)} \bigg( \esup_{t \in [y_k,y_{k+1})} v(t)^{-1} \bigg)^{r'} \bigg)^{ 1 / {r'}} \\
& \hspace{-5cm} \lesssim \sup_{\{b_k\}} \bigg\| \bigg\{ \bigg( \int_{b_{k-1}}^{b_k} \bigg( \int_x^{b_k} w \bigg)^{{r} / {q}} u(x)\,dx \bigg)^{1 / r}  \bigg( \esup_{t \in [b_k,b_{k+1})} v(t)^{-1} \bigg) \bigg\} \bigg\|_{\ell^{\rho} (\Z)} = S.
\end{align*}

Combining, we arrive at
\begin{align*}
\bigg\| \bigg\{ 2^{{k} / {r}} \bigg( \int_{x_k}^{x_{k + 1}} w \bigg) \bigg( \int_{x_k}^{\infty} w(x) \bigg( \esup_{t \in [x,\infty)} v(t)^{-1}\bigg)^{q'}\,dx \bigg)^{1 / {q'}} \bigg\} \bigg\|_{\ell^{\rho}(\zstroke \backslash \{M\})} \le S.
\end{align*}
	
{\rm (ii)} Let $1 \le r$. Let $i \in \Z$. We set 
\begin{equation*}
{\mathcal N}_i^3 : = \left\{ m \in \zstroke:\, 2^i \le \sup_{k < m} 2^{k / r} \bigg( \int_{x_k}^{x_{k+1}} w \bigg) < 2^{i + 1} \right\}
\end{equation*}
and define ${\mathcal N}^3 : = \bigcup_n {\mathcal N}_{i_n}^3 \subset \Z$, where ${\mathcal N}_{i_n}^3 \neq \emptyset$ and ${\mathcal N}_i = \emptyset$ for $i \not \in \{i_n\}$. Next, we define
$$
m_{i_n} = \inf {\mathcal N}_{i_n}^3, \qquad m_{i_n}^+ = \sup {\mathcal N}_{i_n}^3, \qquad z_n = x_{m_{i_n}}.
$$
Suppose, for simplicity that $\card \bigcup_n \{ i_n \} = \aleph_0$.

Applying \cite[Lemma 2.3]{mus.2017}, we get that
\begin{align*}
\bigg\| \bigg\{ 2^{{k} / {r}} \bigg( \int_{x_k}^{x_{k + 1}} w \bigg) \bigg( \int_{x_k}^{\infty} w(x) \bigg( \esup_{t \in [x,\infty)} v(y)^{-1}\bigg)^{q'}\,dx \bigg)^{1 / {q'}} \bigg\} \bigg\|_{\ell^{\rho}(\zstroke \backslash \{M\})} & \\
& \hspace{-7cm} = \sup_{k < M} 2^{k / r} \bigg( \int_{x_k}^{x_{k + 1}} w \bigg) \bigg( \int_{x_k}^{\infty} w(x) \bigg( \esup_{t \in [x,\infty)} v(t)^{-1} \bigg)^{q'} \,dx \bigg)^{1 / {q'}} \\
& \hspace{-7cm} = \sup_{n} \sup_{m_{i_n} \le k \le m_{i_n}^+} 2^{k / r} \bigg( \int_{x_k}^{x_{k + 1}} w \bigg) \bigg( \int_{x_k}^{\infty} w(x) \bigg( \esup_{t \in [x,\infty)} v(t)^{-1} \bigg)^{q'} \,dx \bigg)^{1 / {q'}} \\
& \hspace{-7cm} \le \sup_{n} \bigg( \sup_{m_{i_n} \le k \le m_{i_n}^+} 2^{k / r} \bigg( \int_{x_k}^{x_{k + 1}} w \bigg) \bigg) \bigg( \int_{z_n}^{\infty} w(x) \bigg( \esup_{t \in [x,\infty)} v(t)^{-1} \bigg)^{q'} \,dx \bigg)^{1 / {q'}} \\
& \hspace{-7cm} \lesssim \sup_{n} 2^{i_n} \bigg( \int_{z_n}^{\infty} w(x) \bigg( \esup_{t \in [x,\infty)} v(t)^{-1} \bigg)^{q'} \,dx \bigg)^{1 / {q'}} \\
& \hspace{-7cm} = \sup_{n} 2^{i_n} \bigg( \sum_{ n \le m} \int_{z_m}^{z_{m+1}} w(x) \bigg( \esup_{t \in [x,\infty)} v(t)^{-1} \bigg)^{q'} \,dx \bigg)^{1 / {q'}} \\
& \hspace{-7cm} \approx \sup_{n} 2^{i_n} \bigg( \int_{z_n}^{z_{n+1}} w(x) \bigg( \esup_{t \in [x,\infty)} v(t)^{-1} \bigg)^{q'} \,dx \bigg)^{1 / {q'}} \\
& \hspace{-7cm} \le \sup_{n} 2^{i_n} \bigg( \int_{z_n}^{z_{n+1}} w \bigg)^{1 / {q'}} \bigg( \esup_{t \in [z_n,\infty)} v(t)^{-1} \bigg).
\end{align*}

Let $j \in \Z$. We set 
\begin{equation*}
{\mathcal N}_j^4 : = \left\{ l \in \Z:\, 2^j \le \sup_{n < l} 2^{i_n} \bigg( \int_{z_n}^{z_{n+1}} w \bigg)^{1 / {q'}} < 2^{j + 1} \right\}
\end{equation*}
and define ${\mathcal N}^4 = \bigcup_k {\mathcal N}_{j_k}^4 \subset \Z$, where ${\mathcal N}_{j_k}^4 \neq \emptyset$ and ${\mathcal N}_j^4 = \emptyset$ for $j \not \in \{j_k\}$. Next we define
$$
l_{j_k} = \inf {\mathcal N}_{j_k}^4, \qquad l_{j_k}^+ = \sup {\mathcal N}_{j_k}^4, \qquad y_k = z_{l_{j_k}}.
$$
Assume, for simplicity that $\card \bigcup_k \{ j_k \} = \aleph_0$.

We have, by \cite[Lemma 2.3]{mus.2017}, that
\begin{align*}
\sup_{n} 2^{i_n} \bigg( \int_{z_n}^{z_{n+1}} w \bigg)^{1 / {q'}} \bigg( \esup_{t \in [z_n,\infty)} v(t)^{-1} \bigg) & \\
& \hspace{-5cm} = \sup_k \sup_{l_{j_k} \le n \le l_{j_k}^+} 2^{i_n} \bigg( \int_{z_n}^{z_{n+1}} w \bigg)^{1 / {q'}} \bigg( \esup_{t \in [z_n,\infty)} v(t)^{-1} \bigg) \\
& \hspace{-5cm} \le \sup_k \bigg( \sup_{l_{j_k} \le n \le l_{j_k}^+} 2^{i_n} \bigg( \int_{z_n}^{z_{n+1}} w \bigg)^{1 / {q'}} \bigg) \bigg( \esup_{t \in [z_{l_{j_k}},\infty)} v(t)^{-1} \bigg) \\
& \hspace{-5cm} \lesssim \sup_k 2^{j_k} \bigg( \esup_{t \in [y_k,\infty)} v(t)^{-1} \bigg) \\
& \hspace{-5cm} = \sup_k 2^{j_k} \bigg( \sup_{k \le m} \esup_{t \in [y_m,y_{m+1})} v(t)^{-1} \bigg) \\
& \hspace{-5cm} \approx \sup_k 2^{j_k} \bigg( \esup_{t \in [y_k,y_{k+1})} v(t)^{-1} \bigg).
\end{align*}

Note that $2^{j_k} \le 2 \, \sup_{l_{j_{k-2}} \le n < l_{j_k}} 2^{i_n} \bigg( \int_{z_n}^{z_{n + 1}} w \bigg)^{1 / {q'}}$. Hence
\begin{align*}
\sup_{n} 2^{i_n} \bigg( \int_{z_n}^{z_{n+1}} w \bigg)^{1 / {q'}} \bigg( \esup_{t \in [z_n,\infty)} v(t)^{-1} \bigg) & \\
& \hspace{-5cm} \lesssim \sup_k \bigg( \sup_{l_{j_{k-2}} \le n < l_{j_k}} 2^{i_n} \bigg( \int_{z_n}^{z_{n + 1}} w \bigg)^{1 / {q'}} \bigg)\bigg( \esup_{t \in [y_k,y_{k+1})} v(t)^{-1} \bigg).
\end{align*}

Since $2^{i_n} \le 2\,\sup_{m_{i_{n-2}} \le k < m_{i_n}} 2^{k / r} \bigg( \int_{x_k}^{x_{k+1}} w \bigg)$, we have that
\begin{align*}
2^{i_n} & \lesssim  \sup_{m_{i_{n-2}}\le k < m_{i_n}} \bigg( \int_{x_{k - 1}}^{x_k} u (t)\,dt \bigg)^{1 / r} \bigg( \int_{x_k}^{x_{k+1}} w \bigg) \\
& \le \sup_{m_{i_{n-2}} \le k < m_{i_n}} \bigg( \int_{x_{k - 1}}^{x_k} u (t)\bigg( \int_t^{x_{k+1}} w \bigg)^r \,dt \bigg)^{1 / r} \\
& \le \bigg( \int_{z_{n - 3}}^{z_n} u (t)\bigg( \int_t^{z_n} w \bigg)^r \,dt \bigg)^{1 / r}.
\end{align*}

Hence
\begin{align*}
\sup_{n} 2^{i_n} \bigg( \int_{z_n}^{z_{n+1}} w \bigg)^{1 / {q'}} \bigg( \esup_{t \in [z_n,\infty)} v(t)^{-1} \bigg) & \\
& \hspace{-5cm} \lesssim \sup_k \bigg( \sup_{l_{j_{k-2}} \le n < l_{j_k}} \bigg( \int_{z_{n - 3}}^{z_n} u (t) \bigg( \int_t^{z_n} w \bigg)^r \,dt \bigg)^{1 / r} \bigg( \int_{z_n}^{z_{n + 1}} w \bigg)^{1 / {q'}} \bigg) \bigg( \esup_{t \in [y_k,y_{k+1})} v(t)^{-1} \bigg) \\ 
& \hspace{-5cm} \le \sup_k \bigg( \sup_{l_{j_{k-2}} \le n < l_{j_k}} \bigg( \int_{z_{n - 3}}^{z_{n+1}} u (t)\bigg( \int_t^{z_{n+1}} w \bigg)^r \bigg( \int_t^{z_{n + 1}} w \bigg)^{r / {q'}} \,dt \bigg)^{1 / r} \bigg) \bigg( \esup_{t \in [y_k,y_{k+1})} v(t)^{-1} \bigg) \\
& \hspace{-5cm} = \sup_k \bigg( \sup_{l_{j_{k-2}} \le n < l_{j_k}} \bigg( \int_{z_{n - 3}}^{z_{n+1}} u (t)\bigg( \int_t^{z_{n+1}} w \bigg)^{r / q} \,dt \bigg)^{1 / r} \bigg) \bigg( \esup_{t \in [y_k,y_{k+1})} v(t)^{-1} \bigg) \\
& \hspace{-5cm} \le \sup_k \bigg( \int_{y_{k - 8}}^{y_k} u (t) \bigg( \int_t^{y_k} w \bigg)^{r / q} \,dt \bigg)^{1 / r} \bigg( \esup_{t \in [y_k,y_{k+1})} v(t)^{-1} \bigg) \\
& \hspace{-5cm} \lesssim \sup_{\{b_k\}} \bigg\| \bigg\{ \bigg( \int_{b_{k-1}}^{b_k} \bigg( \int_x^{b_k} w \bigg)^{{r} / {q}} u(x)\,dx \bigg)^{1 / r}  \bigg( \esup_{t \in [b_k,b_{k+1})} v(t)^{-1} \bigg) \bigg\} \bigg\|_{\ell^{\rho} (\Z)} = S.
\end{align*}
	
Combining, we arrive at
\begin{align*}
\bigg\| \bigg\{ 2^{{k} / {r}} \bigg( \int_{x_k}^{x_{k + 1}} w \bigg) \bigg( \int_{x_k}^{\infty} w(x) \bigg( \esup_{t \in [x,\infty)} v(t)^{-1}\bigg)^{q'}\,dx \bigg)^{1 / {q'}} \bigg\} \bigg\|_{\ell^{\rho}(\zstroke \backslash \{M\})} \le S.
\end{align*}

The proof is completed.
\end{proof}
	
Now we prove that conditions $B_i$, $i=\overline{1,4}$ are necessary for inequality \eqref{main}.
\begin{lem}\label{lem.1.9}
	Let $0 < q,\, r < \infty$ and  $u,\,v,\,w \in {\mathcal W}(0,\infty)$. Assume that inequality \eqref{main} holds.
	
	{\rm (i)} If $1 \le \min\{q,\,r\}$, then $B_1 < \infty$. Moreover, if $C$ is the best constant in \eqref{main}, then $B_1 \lesssim C$.
	
	{\rm (ii)} If $r < 1 \le q$, then $B_2 < \infty$. Moreover, if $C$ is the best constant in \eqref{main}, then $B_2 \lesssim C$.
	
	{\rm (iii)} If $q < 1 \le r$, then $B_3 < \infty$. Moreover, if $C$ is the best constant in \eqref{main}, then $B_3 \lesssim C$.
	
	{\rm (iv)} If $\max\{q,\,r\} < 1$, then $B_4 < \infty$. Moreover, if $C$ is the best constant in \eqref{main}, then $B_4 \lesssim C$.
\end{lem}

\begin{proof}
Assume that inequality \eqref{main} holds. Suppose that $\{x_k\}_{k \in \zstroke}$ is a covering sequence.

{\rm (i)} Let $1 \le \min\{q,\,r\}$. In view of Remark \ref{cor.2.0}, $B_1 = D_2$. So, the statement follows by Theorems \ref{thm.1.5} and \ref{cor.1.2}.

{\rm (ii)} Let $r < 1 \le q$. In view of Remark \ref{cor.2.0}, $B_2 = E_1$. Hence, the statement follows by Theorems \ref{thm.1.5} and \ref{cor.1.2}.

{\rm (iii)} Let $q < 1 \le r$. Applying \cite[Lemma 2.3]{mus.2017}, we get that
\begin{align*}
B_3 = & \sup_{k \in \zstroke} \sup_{x_k \le t < x_{k+1}} \bigg( \int_0^t u \bigg)^{1 / r} \bigg( \int_t^{\infty} \bigg( \int_t^x w \bigg)^{q'} w(x) \bigg( \esup_{y \in [x,\infty)} v(y)^{-1} \bigg)^{q'} \,dx\bigg)^{ 1 / {q'}} \\
\lesssim & \sup_{k \in \zstroke} 2^{k / r} \bigg( \int_{x_k}^{\infty} \bigg( \int_{x_k}^x w \bigg)^{q'} w(x) \bigg( \esup_{y \in [x,\infty)} v(y)^{-1} \bigg)^{q'} \,dx\bigg)^{ 1 / {q'}} \\
= & \sup_{k \in \zstroke} 2^{k / r} \bigg( \sum_{m = k}^M \int_{x_m}^{x_{m+1}} \bigg( \int_{x_k}^x w \bigg)^{q'} w(x) \bigg( \esup_{y \in [x,\infty)} v(y)^{-1} \bigg)^{q'} \,dx\bigg)^{ 1 / {q'}} \\
\approx & \sup_{k \in \zstroke} 2^{k / r} \bigg( \sum_{m = k}^M \int_{x_m}^{x_{m+1}} \bigg( \int_{x_m}^x w \bigg)^{q'} w(x) \bigg( \esup_{y \in [x,\infty)} v(y)^{-1} \bigg)^{q'} \,dx\bigg)^{ 1 / {q'}} \\
& +  \sup_{k < M} 2^{k / r} \bigg( \sum_{m = k + 1}^M \bigg( \int_{x_k}^{x_m} w \bigg)^{q'} \int_{x_m}^{x_{m+1}}  w(x) \bigg( \esup_{y \in [x,\infty)} v(y)^{-1} \bigg)^{q'} \,dx\bigg)^{ 1 / {q'}} \\
\approx & \sup_{k \in \zstroke} 2^{k / r} \bigg( \int_{x_k}^{x_{k+1}} \bigg( \int_{x_k}^x w \bigg)^{q'} w(x) \bigg( \esup_{y \in [x,\infty)} v(y)^{-1} \bigg)^{q'} \,dx\bigg)^{ 1 / {q'}} \\
& +  \sup_{k < M} 2^{k / r} \bigg( \sum_{m = k + 1}^M \bigg( \sum_{n = k + 1}^m \int_{x_{n-1}}^{x_n} w \bigg)^{q'} \int_{x_m}^{x_{m+1}}  w(x) \bigg( \esup_{y \in [x,\infty)} v(y)^{-1} \bigg)^{q'} \,dx\bigg)^{ 1 / {q'}} \\
\approx & \sup_{k \in \zstroke} 2^{k / r} \bigg( \int_{x_k}^{x_{k+1}} \bigg( \int_{x_k}^x w \bigg)^{q'} w(x) \bigg( \esup_{y \in [x,x_{k + 1})} v(y)^{-1} \bigg)^{q'} \,dx\bigg)^{ 1 / {q'}} \\
& + \sup_{k < M} 2^{k / r} \bigg( \int_{x_k}^{x_{k+1}} \bigg( \int_{x_k}^x w \bigg)^{q'} w(x) \bigg( \esup_{y \in [x_{k + 1},\infty)} v(y)^{-1} \bigg)^{q'} \,dx\bigg)^{ 1 / {q'}} \\
& +  \sup_{k < M} 2^{k / r} \bigg( \sum_{m = k + 1}^M \bigg( \sum_{n = k + 1}^m \int_{x_{n-1}}^{x_n} w \bigg)^{q'} \int_{x_m}^{x_{m+1}}  w(x) \bigg( \esup_{y \in [x,\infty)} v(y)^{-1} \bigg)^{q'} \,dx\bigg)^{ 1 / {q'}} \\
\lesssim & \sup_{k \in \zstroke} 2^{k / r} \bigg( \int_{x_k}^{x_{k+1}} \bigg( \int_{x_k}^x w \bigg)^{q'} w(x) \bigg( \esup_{y \in [x,x_{k + 1})} v(y)^{-1} \bigg)^{q'} \,dx\bigg)^{ 1 / {q'}} \\
& + \sup_{k < M} 2^{k / r} \bigg( \int_{x_k}^{x_{k+1}} w \bigg) \bigg( \int_{x_k}^{x_{k+1}}  w(x) \bigg( \esup_{y \in [x,\infty)} v(y)^{-1} \bigg)^{q'} \,dx\bigg)^{ 1 / {q'}} \\
& +  \sup_{k < M} 2^{k / r} \bigg( \sum_{m = k + 1}^M \bigg( \sum_{n = k + 1}^m \int_{x_{n-1}}^{x_n} w \bigg)^{q'} \int_{x_m}^{x_{m+1}}  w(x) \bigg( \esup_{y \in [x,\infty)} v(y)^{-1} \bigg)^{q'} \,dx\bigg)^{ 1 / {q'}} \\
= & A_2 + \sup_{k < M} 2^{k / r} \bigg( \int_{x_k}^{x_{k+1}} w \bigg) \bigg( \int_{x_k}^{x_{k+1}}  w(x) \bigg( \esup_{y \in [x,\infty)} v(y)^{-1} \bigg)^{q'} \,dx\bigg)^{ 1 / {q'}} \\
& + \sup_{k < M} 2^{k / r} \bigg( \sum_{m = k + 1}^M \bigg( \sum_{n = k + 1}^m \int_{x_{n-1}}^{x_n} w \bigg)^{q'} \int_{x_m}^{x_{m+1}}  w(x) \bigg( \esup_{y \in [x,\infty)} v(y)^{-1} \bigg)^{q'} \,dx\bigg)^{ 1 / {q'}}.
\end{align*}

Denote by
$$
H_1 : = \sup_{k < M} 2^{k / r} \bigg( \sum_{m = k + 1}^M \bigg( \sum_{n = k + 1}^m \int_{x_{n-1}}^{x_n} w \bigg)^{q'} \int_{x_m}^{x_{m+1}}  w(x) \bigg( \esup_{y \in [x,\infty)} v(y)^{-1} \bigg)^{q'} \,dx\bigg)^{ 1 / {q'}}.
$$

Let $q < 1/2$, which means that $q' < 1$. Using Jensen's inequality, we obtain for any $k < M$ that
\begin{align}
\sum_{m = k + 1}^M \bigg( \sum_{n = k + 1}^m \int_{x_{n-1}}^{x_n} w \bigg)^{q'} \int_{x_m}^{x_{m+1}} w(x) \bigg( \esup_{t \in [x,\infty)} v(t)^{-1} \bigg)^{q'} \,dx & \notag \\
& \hspace{-7cm} \le \sum_{m = k + 1}^M \sum_{n = k + 1}^m \bigg( \int_{x_{n-1}}^{x_n} w \bigg)^{q'} \int_{x_m}^{x_{m+1}} w(x) \bigg( \esup_{t \in [x,\infty)} v(t)^{-1} \bigg)^{q'} \,dx \notag \\
& \hspace{-7cm} \le \sum_{n = k + 1}^M \sum_{m = n}^M \bigg( \int_{x_{n-1}}^{x_n} w \bigg)^{q'} \int_{x_m}^{x_{m+1}} w(x) \bigg( \esup_{t \in [x,\infty)} v(t)^{-1} \bigg)^{q'} \,dx \notag \\
& \hspace{-7cm} = \sum_{n = k + 1}^M \bigg( \int_{x_{n-1}}^{x_n} w \bigg)^{q'} \sum_{m = n}^M \int_{x_m}^{x_{m+1}} w(x) \bigg( \esup_{t \in [x,\infty)} v(t)^{-1} \bigg)^{q'} \,dx \notag \\
& \hspace{-7cm} = \sum_{n = k + 1}^M \bigg( \int_{x_{n-1}}^{x_n} w \bigg)^{q'} \int_{x_n}^{\infty} w(x) \bigg( \esup_{t \in [x,\infty)} v(t)^{-1} \bigg)^{q'} \,dx. \label{eq.1.1.1}
\end{align}

Thus, by \cite[Lemma 2.3]{mus.2017}, we arrive at
\begin{align*}
H_1  & \le \sup_{k < M} 2^{k / r} \bigg( \sum_{n = k + 1}^M \bigg( \int_{x_{n-1}}^{x_n} w \bigg)^{q'} \int_{x_n}^{\infty} w(x) \bigg( \esup_{t \in [x,\infty)} v(t)^{-1} \bigg)^{q'} \,dx \bigg)^{ 1 / {q'}} \\
& \lesssim \sup_{k < M} 2^{k / r} \bigg( \int_{x_k}^{x_{k + 1}} w \bigg) \bigg( \int_{x_{k+1}}^{\infty} w(x) \bigg( \esup_{t \in [x,\infty)} v(t)^{-1} \bigg)^{q'} \,dx \bigg)^{ 1 / {q'}}.
\end{align*}

Let now $q \ge 1/2$, which means that $q' \ge 1$. Then, by Minkowski's inequality, we have for any $k < M$ that
\begin{align}
\sum_{m = k + 1}^M \bigg( \sum_{n = k + 1}^m \int_{x_{n-1}}^{x_n} w \bigg)^{q'} \int_{x_m}^{x_{m+1}} w(x) \bigg( \esup_{t \in [x,\infty)} v(t)^{-1} \bigg)^{q'} \,dx & \notag \\
& \hspace{-7cm} = \left\{ \left[\sum_{m = k + 1}^M \bigg( \sum_{n = k + 1}^m \bigg( \int_{x_{n-1}}^{x_n} w \bigg) \bigg( \int_{x_m}^{x_{m+1}} w(x) \bigg( \esup_{t \in [x,\infty)} v(t)^{-1} \bigg)^{q'} \,dx \bigg)^{1 / {q'}} \bigg)^{q'} \right]^{1 / {q'}}  \right\}^{q'} \notag \\
& \hspace{-7cm} \le \left\{ \sum_{n = k + 1}^M \left[ \sum_{m = n}^M \bigg( \int_{x_{n-1}}^{x_n} w \bigg)^{q'} \bigg( \int_{x_m}^{x_{m+1}} w(x) \bigg( \esup_{t \in [x,\infty)} v(t)^{-1} \bigg)^{q'} \,dx \bigg)  \right]^{1 / {q'}}  \right\}^{q'} \notag\\
& \hspace{-7cm} = \left\{ \sum_{n = k + 1}^M \bigg( \int_{x_{n-1}}^{x_n} w \bigg) \left[ \sum_{m = n}^M \bigg( \int_{x_m}^{x_{m+1}} w(x) \bigg( \esup_{t \in [x,\infty)} v(t)^{-1} \bigg)^{q'} \,dx \bigg)  \right]^{1 / {q'}}  \right\}^{q'} \notag \\
& \hspace{-7cm} = \left\{ \sum_{n = k + 1}^M \bigg( \int_{x_{n-1}}^{x_n} w \bigg) \bigg( \int_{x_{n}}^{\infty} w(x) \bigg( \esup_{t \in [x,\infty)} v(t)^{-1} \bigg)^{q'} \,dx \bigg)^{1 / {q'}}  \right\}^{q'}. \label{eq.1.1.2}
\end{align}

Hence, by \cite[Lemma 2.3]{mus.2017}, we get that
\begin{align*}
H_1 & \le \sup_{k < M} 2^{k / r} \sum_{n = k + 1}^M \bigg( \int_{x_{n-1}}^{x_n} w \bigg) \bigg( \int_{x_{n}}^{\infty} w(x) \bigg( \esup_{t \in [x,\infty)} v(t)^{-1} \bigg)^{q'} \,dx \bigg)^{1 / {q'}} \\
& \approx \sup_{k < M} 2^{k / r} \bigg( \int_{x_k}^{x_{k + 1}} w \bigg) \bigg( \int_{x_{k+1}}^{\infty} w(x) \bigg( \esup_{t \in [x,\infty)} v(t)^{-1} \bigg)^{q'} \,dx \bigg)^{1 / {q'}}.
\end{align*}

Thus, in both cases, we have that 
\begin{align*}
H_1 & \lesssim \sup_{k < M} 2^{k / r} \bigg( \int_{x_k}^{x_{k + 1}} w \bigg) \bigg( \int_{x_{k+1}}^{\infty} w(x) \bigg( \esup_{t \in [x,\infty)} v(t)^{-1} \bigg)^{q'} \,dx \bigg)^{1 / {q'}}.
\end{align*}

By Lemma \ref{lem.2.5}, we obtain that
\begin{align*}
\sup_{k < M} 2^{k / r} \bigg( \int_{x_k}^{x_{k+1}} w \bigg) \bigg( \int_{x_k}^{x_{k+1}}  w(x) \bigg( \esup_{y \in [x,\infty)} v(y)^{-1} \bigg)^{q'} \,dx\bigg)^{ 1 / {q'}} & \\
& \hspace{-7.5cm} + \sup_{k < M} 2^{k / r} \bigg( \sum_{m = k + 1}^M \bigg( \sum_{n = k + 1}^m \int_{x_{n-1}}^{x_n} w \bigg)^{q'} \int_{x_m}^{x_{m+1}}  w(x) \bigg( \esup_{y \in [x,\infty)} v(y)^{-1} \bigg)^{q'} \,dx\bigg)^{ 1 / {q'}} \\
& \hspace{-8cm} = \sup_{k < M} 2^{k / r} \bigg( \int_{x_k}^{x_{k+1}} w \bigg) \bigg( \int_{x_k}^{x_{k+1}}  w(x) \bigg( \esup_{y \in [x,\infty)} v(y)^{-1} \bigg)^{q'} \,dx\bigg)^{ 1 / {q'}} + H_1 \\
& \hspace{-8cm} \lesssim \sup_{k < M} 2^{k / r} \bigg( \int_{x_k}^{x_{k+1}} w \bigg) \bigg( \int_{x_k}^{x_{k+1}}  w(x) \bigg( \esup_{y \in [x,\infty)} v(y)^{-1} \bigg)^{q'} \,dx\bigg)^{ 1 / {q'}} \\
& \hspace{-7.5cm} + \sup_{k < M} 2^{k / r} \bigg( \int_{x_k}^{x_{k + 1}} w \bigg) \bigg( \int_{x_{k+1}}^{\infty} w(x) \bigg( \esup_{t \in [x,\infty)} v(t)^{-1} \bigg)^{q'} \,dx \bigg)^{1 / {q'}} \\
& \hspace{-8cm} \lesssim \sup_{k < M} 2^{k / r} \bigg( \int_{x_k}^{x_{k + 1}} w \bigg) \bigg( \int_{x_k}^{\infty} w(x) \bigg( \esup_{t \in [x,\infty)} v(t)^{-1} \bigg)^{q'} \,dx \bigg)^{1 / {q'}} \lesssim S.
\end{align*}

Combining,  we arrive at
$$
B_3 \lesssim A_2 + S,
$$
and the statement follows by Theorem \ref{thm.1.5}, Lemma \ref{lem.1.6} and Lemma \ref{lem.2.3}.

{\rm (iv)} Let $\max\{q,\,r\} < 1$. Applying \cite[Lemma 2.3]{mus.2017}, we get that
\begin{align*}
B_4 = & \bigg( \sum_{k \in \zstroke} \int_{x_k}^{x_{k+1}} \bigg( \int_0^t u \bigg)^{r'} u(t) \bigg( \int_t^{\infty} \bigg( \int_t^x w \bigg)^{q'} w(x) \bigg( \esup_{t \in [x,\infty)} v(t)^{-1} \bigg)^{q'} \,dx\bigg)^{ {r'} / {q'}}\,dt \bigg)^{ 1 / {r'}} \\
\le & \bigg( \sum_{k \in \zstroke} \int_{x_k}^{x_{k+1}} \bigg( \int_0^t u \bigg)^{r'} u(t) \,dt \cdot \bigg( \int_{x_k}^{\infty} \bigg( \int_{x_k}^x w \bigg)^{q'} w(x) \bigg( \esup_{t \in [x,\infty)} v(t)^{-1} \bigg)^{q'} \,dx\bigg)^{ {r'} / {q'}} \bigg)^{ 1 / {r'}} \\
\approx & \bigg( \sum_{k \in \zstroke} 2^{k(r' + 1)} \bigg( \int_{x_k}^{\infty} \bigg( \int_{x_k}^x w \bigg)^{q'} w(x) \bigg( \esup_{t \in [x,\infty)} v(t)^{-1} \bigg)^{q'} \,dx\bigg)^{ {r'} / {q'}} \bigg)^{ 1 / {r'}} \\
= & \bigg( \sum_{k \in \zstroke} 2^{k(r' + 1)} \bigg( \sum_{m = k}^M \int_{x_m}^{x_{m+1}} \bigg( \int_{x_k}^x w \bigg)^{q'} w(x) \bigg( \esup_{t \in [x,\infty)} v(t)^{-1} \bigg)^{q'} \,dx\bigg)^{ {r'} / {q'}} \bigg)^{ 1 / {r'}} \\
\approx & \bigg( \sum_{k \in \zstroke} 2^{k(r' + 1)} \bigg( \sum_{m = k}^M \int_{x_m}^{x_{m+1}} \bigg( \int_{x_m}^x w \bigg)^{q'} w(x) \bigg( \esup_{t \in [x,\infty)} v(t)^{-1} \bigg)^{q'} \,dx\bigg)^{ {r'} / {q'}} \bigg)^{ 1 / {r'}} \\
& + \bigg( \sum_{k < M} 2^{k(r' + 1)} \bigg( \sum_{m = k + 1}^M \bigg( \int_{x_k}^{x_m} w \bigg)^{q'} \int_{x_m}^{x_{m+1}} w(x) \bigg( \esup_{t \in [x,\infty)} v(t)^{-1} \bigg)^{q'} \,dx\bigg)^{ {r'} / {q'}} \bigg)^{ 1 / {r'}} \\
\approx & \bigg( \sum_{k \in \zstroke} 2^{k(r' + 1)} \bigg( \int_{x_k}^{x_{k+1}} \bigg( \int_{x_k}^x w \bigg)^{q'} w(x) \bigg( \esup_{t \in [x,\infty)} v(t)^{-1} \bigg)^{q'} \,dx\bigg)^{ {r'} / {q'}} \bigg)^{ 1 / {r'}} \\
& + \bigg( \sum_{k < M} 2^{k(r' + 1)} \bigg( \sum_{m = k + 1}^M \bigg( \sum_{n = k + 1}^m \int_{x_{n-1}}^{x_n} w \bigg)^{q'} \int_{x_m}^{x_{m+1}} w(x) \bigg( \esup_{t \in [x,\infty)} v(t)^{-1} \bigg)^{q'} \,dx\bigg)^{ {r'} / {q'}} \bigg)^{ 1 / {r'}} \\
\approx & \bigg( \sum_{k \in \zstroke} 2^{k(r' + 1)} \bigg( \int_{x_k}^{x_{k+1}} \bigg( \int_{x_k}^x w \bigg)^{q'} w(x) \bigg( \esup_{t \in [x,x_{k+1})} v(t)^{-1} \bigg)^{q'} \,dx\bigg)^{ {r'} / {q'}} \bigg)^{ 1 / {r'}} \\
& + \bigg( \sum_{k < M} 2^{k(r' + 1)} \bigg( \int_{x_k}^{x_{k+1}} \bigg( \int_{x_k}^x w \bigg)^{q'} w(x) \bigg( \esup_{t \in [x_{k+1},\infty)} v(t)^{-1} \bigg)^{q'} \,dx \bigg)^{ {r'} / {q'}} \bigg)^{ 1 / {r'}} \\
& + \bigg( \sum_{k < M} 2^{k(r' + 1)} \bigg( \sum_{m = k + 1}^M \bigg( \sum_{n = k + 1}^m \int_{x_{n-1}}^{x_n} w \bigg)^{q'} \int_{x_m}^{x_{m+1}} w(x) \bigg( \esup_{t \in [x,\infty)} v(t)^{-1} \bigg)^{q'} \,dx\bigg)^{ {r'} / {q'}} \bigg)^{ 1 / {r'}} \\
\le & \bigg( \sum_{k \in \zstroke} 2^{k(r' + 1)} \bigg( \int_{x_k}^{x_{k+1}} \bigg( \int_{x_k}^x w \bigg)^{q'} w(x) \bigg( \esup_{t \in [x,x_{k+1})} v(t)^{-1} \bigg)^{q'} \,dx\bigg)^{ {r'} / {q'}} \bigg)^{ 1 / {r'}} \\
& + \bigg( \sum_{k < M} 2^{k(r' + 1)} \bigg( \bigg( \int_{x_k}^{x_{k+1}} w \bigg)^{q'} \bigg( \int_{x_k}^{x_{k+1}}  w(x) \bigg( \esup_{t \in [x,\infty)} v(t)^{-1} \bigg)^{q'} \,dx \bigg) \bigg)^{ {r'} / {q'}} \bigg)^{ 1 / {r'}} \\
& + \bigg( \sum_{k < M} 2^{k(r' + 1)} \bigg( \sum_{m = k + 1}^M \bigg( \sum_{n = k + 1}^m \int_{x_{n-1}}^{x_n} w \bigg)^{q'} \int_{x_m}^{x_{m+1}} w(x) \bigg( \esup_{t \in [x,\infty)} v(t)^{-1} \bigg)^{q'} \,dx\bigg)^{ {r'} / {q'}} \bigg)^{ 1 / {r'}} \\
= & A_2 + \bigg( \sum_{k < M} 2^{k(r' + 1)} \bigg( \int_{x_k}^{x_{k+1}} w \bigg)^{r'} \bigg( \int_{x_k}^{x_{k+1}}  w(x) \bigg( \esup_{t \in [x,\infty)} v(t)^{-1} \bigg)^{q'} \,dx \bigg)^{ {r'} / {q'}} \bigg)^{ 1 / {r'}} \\
& + \bigg( \sum_{k < M} 2^{k(r' + 1)} \bigg( \sum_{m = k + 1}^M \bigg( \sum_{n = k + 1}^m \int_{x_{n-1}}^{x_n} w \bigg)^{q'} \int_{x_m}^{x_{m+1}} w(x) \bigg( \esup_{t \in [x,\infty)} v(t)^{-1} \bigg)^{q'} \,dx \bigg)^{ {r'} / {q'}} \bigg)^{ 1 / {r'}}.
\end{align*}

Denote by
$$
H_2 : = \bigg( \sum_{k < M} 2^{k(r' + 1)} \bigg( \sum_{m = k + 1}^M \bigg( \sum_{n = k + 1}^m \int_{x_{n-1}}^{x_n} w \bigg)^{q'} \int_{x_m}^{x_{m+1}} w(x) \bigg( \esup_{t \in [x,\infty)} v(t)^{-1} \bigg)^{q'} \,dx\bigg)^{ {r'} / {q'}} \bigg)^{ 1 / {r'}}.
$$

Let $q < 1/2$. Applying inequality \eqref{eq.1.1.1}, by \cite[Lemma 2.3]{mus.2017}, we arrive at
\begin{align*}
H_2  & \le \bigg( \sum_{k < M} 2^{k(r' + 1)} \bigg( \sum_{n = k + 1}^M \bigg( \int_{x_{n-1}}^{x_n} w \bigg)^{q'} \int_{x_n}^{\infty} w(x) \bigg( \esup_{t \in [x,\infty)} v(t)^{-1} \bigg)^{q'} \,dx \bigg)^{ {r'} / {q'}} \bigg)^{ 1 / {r'}} \\
& \lesssim \bigg( \sum_{k < M} 2^{k(r' + 1)} \bigg( \int_{x_k}^{x_{k+1}} w \bigg)^{r'} \bigg( \int_{x_{k+1}}^{\infty} w(x) \bigg( \esup_{t \in [x,\infty)} v(t)^{-1} \bigg)^{q'} \,dx \bigg)^{ {r'} / {q'}} \bigg)^{ 1 / {r'}}.
\end{align*}

Let now $q \ge 1/2$. Applying \eqref{eq.1.1.2}, by \cite[Lemma 2.3]{mus.2017}, we get that
\begin{align*}
H_2  & \le \bigg( \sum_{k < M} 2^{k(r' + 1)} \bigg( \sum_{n = k + 1}^M \bigg( \int_{x_{n-1}}^{x_n} w \bigg) \bigg( \int_{x_{n}}^{\infty} w(x) \bigg( \esup_{t \in [x,\infty)} v(t)^{-1} \bigg)^{q'} \,dx \bigg)^{1 / {q'}}\bigg)^{ {r'}} \bigg)^{ 1 / {r'}} \\
& \lesssim \bigg( \sum_{k < M} 2^{k(r' + 1)} \bigg( \int_{x_k}^{x_{k+1}} w \bigg)^{r'} \bigg( \int_{x_{k+1}}^{\infty} w(x) \bigg( \esup_{t \in [x,\infty)} v(t)^{-1} \bigg)^{q'} \,dx \bigg)^{ {r'} / {q'}} \bigg)^{ 1 / {r'}}.
\end{align*}

Thus, in both cases, we have that 
\begin{align*}
H_2  \lesssim \bigg( \sum_{k < M} 2^{k(r' + 1)} \bigg( \int_{x_k}^{x_{k+1}} w \bigg)^{r'} \bigg( \int_{x_{k+1}}^{\infty} w(x) \bigg( \esup_{t \in [x,\infty)} v(t)^{-1} \bigg)^{q'} \,dx \bigg)^{ {r'} / {q'}} \bigg)^{ 1 / {r'}}.
\end{align*}

By Lemma \ref{lem.2.5}, we obtain that
\begin{align*}
\bigg( \sum_{k < M} 2^{k(r' + 1)} \bigg( \int_{x_k}^{x_{k+1}} w \bigg)^{r'} \bigg( \int_{x_k}^{x_{k+1}}  w(x) \bigg( \esup_{t \in [x,\infty)} v(t)^{-1} \bigg)^{q'} \,dx \bigg)^{ {r'} / {q'}} \bigg)^{ 1 / {r'}} & \\
& \hspace{-8.5cm} + \bigg( \sum_{k < M} 2^{k(r' + 1)} \bigg( \sum_{m = k + 1}^M \bigg( \sum_{n = k + 1}^m \int_{x_{n-1}}^{x_n} w \bigg)^{q'} \int_{x_m}^{x_{m+1}} w(x) \bigg( \esup_{t \in [x,\infty)} v(t)^{-1} \bigg)^{q'} \,dx \bigg)^{ {r'} / {q'}} \bigg)^{ 1 / {r'}} \\
& \hspace{-9cm} = \bigg( \sum_{k < M} 2^{k(r' + 1)} \bigg( \int_{x_k}^{x_{k+1}} w \bigg)^{r'} \bigg( \int_{x_k}^{x_{k+1}}  w(x) \bigg( \esup_{t \in [x,\infty)} v(t)^{-1} \bigg)^{q'} \,dx \bigg)^{ {r'} / {q'}} \bigg)^{ 1 / {r'}} + H_2 \\
& \hspace{-9cm} = \bigg( \sum_{k < M} 2^{k(r' + 1)} \bigg( \int_{x_k}^{x_{k+1}} w \bigg)^{r'} \bigg( \int_{x_k}^{x_{k+1}}  w(x) \bigg( \esup_{t \in [x,\infty)} v(t)^{-1} \bigg)^{q'} \,dx \bigg)^{ {r'} / {q'}} \bigg)^{ 1 / {r'}} \\
& \hspace{-8.5cm} + \bigg( \sum_{k < M} 2^{k(r' + 1)} \bigg( \int_{x_k}^{x_{k+1}} w \bigg)^{r'} \bigg( \int_{x_{k+1}}^{\infty} w(x) \bigg( \esup_{t \in [x,\infty)} v(t)^{-1} \bigg)^{q'} \,dx \bigg)^{ {r'} / {q'}} \bigg)^{ 1 / {r'}} \\
& \hspace{-9cm} \lesssim \bigg( \sum_{k < M} 2^{k(r' + 1)} \bigg( \int_{x_k}^{x_{k+1}} w \bigg)^{r'} \bigg( \int_{x_k}^{\infty} w(x) \bigg( \esup_{t \in [x,\infty)} v(t)^{-1} \bigg)^{q'} \,dx \bigg)^{ {r'} / {q'}} \bigg)^{ 1 / {r'}} \lesssim S.
\end{align*}

Combining,  we arrive at
$$
B_4 \lesssim A_2 + S,
$$
and the statement follows by Theorem \ref{thm.1.5}, Lemma \ref{lem.1.6} and Lemma \ref{lem.2.3}.

The proof is completed.
\end{proof}


\section{Proof of the main statement.}\label{s.8}

We are now in a position to prove the main statement.

\noindent{\bf Proof of Theorem \ref{main1}.} 

Note that $F_1 = D$, $F_2 = E_1$ and $F_3 = E_2$.

Since for all $t > 0$
\begin{align*}
\esup_{s \in [t,\infty)} \bigg( \int_t^s w \bigg)^{1 / q} \bigg(\esup_{y\in [s,\infty)} v(y)^{-1} \bigg) & \\
& \hspace{-3cm} \approx \esup_{s \in [t,\infty)} \bigg( \int_t^s \bigg( \int_t^x w \bigg)^{q'} w(x)\,dx \bigg)^{1 / q'} \bigg(\esup_{y\in [s,\infty)} v(y)^{-1} \bigg) \\
& \hspace{-3cm} \le \esup_{s \in [t,\infty)} \bigg( \int_t^s \bigg( \int_t^x w \bigg)^{q'} w(x) \bigg(\esup_{y\in [x,\infty)} v(y)^{-1} \bigg)^{q'} \,dx \bigg)^{1 / q'} \\
& \hspace{-3cm} = \bigg( \int_t^{\infty} \bigg( \int_t^x w \bigg)^{q'} w(x) \bigg(\esup_{y\in [x,\infty)} v(y)^{-1} \bigg)^{q'} \,dx \bigg)^{1 / q'},
\end{align*}
then
\begin{align*}
F_1 = B_1 & = \esup_{t \in (0,\infty)}  \bigg( \int_0^t u \bigg)^{1/r}  \bigg( \esup_{s \in [t,\infty)} \bigg( \int_t^s w \bigg)^{1/q} v(s)^{-1} \bigg) \\
& \le \esup_{t\in (0,\infty)} \bigg(\int_{0}^{t} u \bigg)^{{1} / {r}} \bigg(\int_{t}^{\infty} \bigg(\int_{t}^{x} w \bigg)^{q'} w(x) \bigg(\esup_{y\in (x,\infty)} v(y)^{-1} \bigg)^{q'} \, dx \bigg)^{1 / q'} = B_3 = F_4
\end{align*}
and
\begin{align*}
F_2 = B_2 & = \bigg( \int_0^{\infty}  \bigg( \int_0^t u \bigg)^{r'}  u(t) \bigg( \esup_{s \in [t,\infty)} \bigg( \int_t^s w \bigg)^{1 / q} v(s)^{-1} \bigg)^{r'} \,dt\bigg)^{1/r'} \\
& \le \bigg( \int_0^{\infty}  \bigg( \int_0^t u \bigg)^{r'}  u(t) \bigg( \int_t^{\infty} \bigg( \int_t^x w \bigg)^{q'} w(x) \bigg(\esup_{y\in (x,\infty)} v(y)^{-1} \bigg)^{q'}\,dx \bigg)^{r' / q'} \,dt\bigg)^{1/r'} = B_4 = F_5.
\end{align*}

So, the statement follows from Theorem \ref{thm.1.5}, Theorem \ref{cor.1.2}, Lemma \ref{lem.1.6}, Lemma \ref{lem.1.8} and Lemma \ref{lem.1.9}. \qed


\end{document}